\documentclass[10pt]{amsart}

\usepackage{amssymb}
\usepackage{amscd}
\usepackage{upref}

\theoremstyle{plain}
\newtheorem{thm}{Theorem}[section]

\newtheorem{lem}[thm]{Lemma}
\newtheorem{prop}[thm]{Proposition}

\theoremstyle{definition}
\newtheorem{defn}[thm]{Definition}

\theoremstyle{remark}

\numberwithin{equation}{section}

\begin{document}
\title[A Proof of BV Identity ]
{A Homotopy Theoretic Proof of the BV Identity in Loop Homology} 
\author{Hirotaka Tamanoi} 
\address[] {Department of Mathematics, 
University of California Santa Cruz \newline 
\indent Santa Cruz, CA 95064} 
\email[]{tamanoi@math.ucsc.edu} 
\date{}
\begin{abstract}
Chas and Sullivan proved the existence of a Batalin-Vilkovisky algebra
structure in the homology of free loop spaces on closed finite
dimensional smooth manifolds using chains and chain homotopies. This
algebraic structure involves an associative product called the loop
product, a Lie bracket called the loop bracket, and a square 0
operator called the BV operator. Cohen and Jones gave a homotopy
theoretic description of the loop product in terms of spectra. In this
paper, we give an explicit homotopy theoretic description of the loop bracket
and, using this description, we give a homological proof of the BV
identity connecting the loop product, the loop bracket, and the BV
operator. The proof is based on an observation that the loop bracket
and the BV derivation are given by the same cycle in the free loop
space, except that they differ by parametrization of loops. 
\end{abstract}
\maketitle

\tableofcontents

\section{Introduction}

Let $M$ be a closed oriented smooth $d$-manifold, and let $LM$ be its free
loop space of continuous maps from $S^1=\mathbb{R}/\mathbb{Z}$ to
$M$. Chas and Sullivan \cite{CS} proved that its homology $H_*(LM)$
has the structure of a Batalin-Vilkovisky (BV) algebra. Namely, they
showed that the loop homology has an associative graded commutative
product $\cdot$ of degree $-d$ called the loop product, a compatible
Lie bracket $\{\ ,\ \}$ of degree $1$ called the loop bracket, and an
operator $\Delta$ of degree $1$ coming from the circle action on $LM$
called the BV operator, satisfying the following BV identity for
$a,b\in H_*(LM)$:
\begin{equation}\label{bv identity}
\Delta(a\cdot b)=(\Delta a)\cdot b+(-1)^{|a|-d}a\cdot\Delta b
+(-1)^{|a|-d}\{a,b\},
\end{equation}
where $|a|$ denotes the homological degree in $H_*(LM)$. They 
showed the existence of a BV algebra structure using transversal
chains and chain homotopies. 

Cohen and Jones \cite{CJ} gave a homotopy theoretic description of the
loop product using Pontrjagin-Thom construction and showed that the
spectrum $LM^{-TM}$ is a ring spectrum with respect to the loop
product. Voronov \cite{V} showed that the homology of the cacti
operad acts on the loop homology, which automatically implies that the
loop homology has a BV algebra structure by a general theorem due to
\cite{G}. In this context, cycles in the cacti operad give rise to
homology operator on $H_*(LM)$, and homologous relations among cycles
give rise to identities satisfied by corresponding homology operators. Chas
and Sullivan explicitly constructed these cycles and homologous
relations among them on chain level. The above BV identity was proved
in this way.

In this paper, after reviewing homotopy theoretic description of the
loop product due to \cite{CJ}, we give explicit homotopy theoretic
reformulation of the loop bracket whose chain description was given in
\cite{CS}. We then give homological proof of the BV identity using
this description. Our main observation is that cycles representing
$(\Delta a)\cdot b+(-1)^{|a|-d}a\cdot\Delta b$ and
$(-1)^{|a|-d}\{a,b\}$ can be taken to be the same up to
reparametrization of loops, and the correction of this difference of
parametrization of loops yields the term $\Delta(a\cdot b)$, proving
the BV identity. Here, we give an outline of the proof. Details are given in subsequent sections. 

\begin{proof}[Outline of the homotopy theoretic proof of BV identity \eqref{bv identity}]
Let $\phi: M \rightarrow M\times M$ be the diagonal map, and let $e: S^1\times LM\times LM \longrightarrow M\times M$ be an evaluation map 
\begin{equation}
e(t,\gamma,\eta)=
\begin{cases}
\bigl(\gamma(0),\eta(2t)\bigr)   & 0\le t\le \tfrac12, \\
\bigl(\gamma(2t-1),\eta(1)\bigr) & \tfrac12\le t\le 1. 
\end{cases}
\end{equation}
Let $j: e^{-1}\bigl(\phi(M)\bigr) \longrightarrow S^1\times LM\times
LM$ be the inclusion map. Then there exist maps $\iota_1, \iota_2:
e^{-1}\bigl(\phi(M)\bigr) \longrightarrow LM$ (see \eqref{iota1-1},
\eqref{iota1-2}, \eqref{iota2-1}, \eqref{iota2-2}) such that
\begin{align} 
{\iota_1}_*j_!([S^1]\times a\times b)&=(-1)^{|a|+d(|a|-d)}\{a,b\}, \\
{\iota_2}_*j_!([S^1]\times a\times b)&=(-1)^{|a|+d(|a|-d)}a\cdot\Delta
b + (-1)^{d|a|}(\Delta a)\cdot b,
\end{align}
where $a,b\in H_*(LM)$. Note that the loop bracket and the BV 
derivation are defined on the same subset $e^{-1}\bigl(\phi(M)\bigr)$, 
with different interaction maps $\iota_1$ and $\iota_2$.  
Here, the loop $\iota_2(t,\gamma,\eta)$
coincides with $\iota_1(t,\gamma,\eta)$ rotated by $t$ (Lemma
\ref{rotation}). Adjusting the difference of parametrization yields
the term $\Delta(a\cdot b)$ (Proposition \ref{Pi Iota1} and Theorem \ref{proof of BV}), and completes the proof of BV identity.
\end{proof} 

The organization of this paper is as follows. After reviewing 
the loop product in \S 2, we give a homotopy theoretic description 
of the loop bracket in \S 3, followed by a homotopy theoretic 
description of the BV derivation in \S 4. In \S 5, we compare 
these two descriptions and prove the BV identity.

\section{The Loop Product}

In this section, we give a quick review of the homotopy theoretic
description of the loop product given in \cite{CJ}. Let $p:
LM\longrightarrow M$ be the base point map given by
$p(\gamma)=\gamma(0)$ for $\gamma\in LM$, and let $\phi: M
\longrightarrow M\times M$ be the diagonal map.  Let $LM\times_M LM$
be the space $(p\times p)^{-1}\bigl(\phi(M)\bigr)$ consisting of pairs
of loops $(\gamma, \eta)$ with the same base points, and let $\iota:
LM\times_M LM \longrightarrow LM$ be the usual loop multiplication map
$\iota(\gamma,\eta)=\gamma\cdot\eta$. Thus we have the following
diagram.
\begin{equation}
\begin{CD}
LM\times LM @<{j}<< LM\times_M LM @>{\iota}>> LM \\
@V{p\times p}VV  @V{q}VV   @.  \\
M\times M @<{\phi}<< M  @.   
\end{CD}
\end{equation}
where $j$ is the inclusion map and $q$ is the restriction of $p\times p$. 

\begin{defn}
The loop product of $a,b\in H_*(LM)$ is given by 
\begin{equation}
a\cdot b=(-1)^{d(|a|-d)}\iota_*j_!(a\times b),
\end{equation}
where $j_!$ is the transfer map of degree $-d$. 
\end{defn}
We recall the construction of the transfer map. This will serve 
as a preparation for a homotopy theoretic description of 
the loop bracket in the next section. Let $\nu$ be the
normal bundle to $\phi(M)$ in $M\times M$. We orient $\nu$ so that we
have $\nu\oplus\phi_*(TM)\cong T(M\times M)|_{\phi(M)}$. Let $u'$ be
the Thom class of $\nu$ with this orientation. Let $N$ be a closed
tubular neighborhood so that we have $D(\nu)\cong N$, where $D(\nu)$
is the closed disc bundle associated to $\nu$. We have $u'\in
\tilde{H}^d(N/\partial N)$. Let $\pi: N\longrightarrow \phi(M)$ be the
projection map, which is also a homotopy equivalence. Let $c:M\times M
\longrightarrow N/\partial N$ be the Thom collapse map, and let
$u=c^*(u')\in H^d(M\times M)$. Then $u$ is the cohomology class dual
to the diagonal in the sense that $u\cap[M\times M]=\phi_*([M])$. Let
$\widetilde{N}=(p\times p)^{-1}(N)$ and let $\tilde{c}: LM\times LM
\longrightarrow \widetilde{N}/\partial\widetilde{N}$ be the Thom collapse
map. The projection map $\pi$ can be lifted to $\tilde{\pi}: \widetilde{N}
\longrightarrow LM\times_M LM$ using the bundle structure of $N$, and
$\tilde{\pi}$ is a homotopy equivalence. Let $\tilde{u}'=(p\times
p)^*(u')\in \tilde{H}^d(\widetilde{N}/\partial\widetilde{N})$ and
$\tilde{u}=(p\times p)^*(u)\in H^d(LM\times LM)$. Now the transfer map
$j_!$ can be defined as the composition of the following maps.

\begin{equation} 
j_!: H_*(LM\times LM) \xrightarrow{\tilde{c}_*}
\tilde{H}_*(\widetilde{N}/\partial\widetilde{N}) \xrightarrow{\tilde{u}'\cap(\
)} H_{*-d}(\widetilde{N}) \xrightarrow[\cong]{\tilde{\pi}_*}
H_{*-d}(LM\underset{M}{\times} LM).
\end{equation}

Geometrically, the map $j_!$ is given by taking transversal
intersection of the cycle representing $a\times b$ with the
codimension $d$ submanifold $LM\times_M LM$ with an appropriate
orientation, and then taking its homology class in $LM\times_M
LM$. The following property of the transfer map is useful.

\begin{prop}\label{jj-formula}
For $a,b\in H_*(LM\times LM)$, we have 
\begin{equation}
j_*j_!(a\times b)=\tilde{u}\cap(a\times b).
\end{equation}
\end{prop}
\begin{proof}
We consider the following commutative diagram with obvious maps. 
\begin{equation*}
\begin{CD}
H^d((LM)^2) @<{\tilde{k}^*}<<  H^d((LM)^2, (LM)^2-LM\underset{M}{\times}LM)
@>{\iota_{\widetilde{N}}^*}>{\cong}>  H^d(\widetilde{N}, \partial\widetilde{N}) \\
@A{(p\times p)^*}AA   @A{(p\times p)^*}AA   @A{(p\times p)^*}AA    \\
H^d(M\times M) @<{k^*}<<   H^d\bigl(M\times M, M\times M-\phi(M)\bigr) 
@>{\iota_N^*}>{\cong}>  H^d(N,\partial N) 
\end{CD}
\end{equation*}
Here $c^*=k^*(\iota_N^*)^{-1}$ and
$\tilde{c}^*=\tilde{k}^*(\iota_{\widetilde{N}}^*)^{-1}$. We let $u''\in
H^d(M\times M, M\times M-\phi(M))$ be the Thom class corresponding to
$u'\in H^d(N,\partial N)$, and let $\tilde{u}''=(p\times
p)^*(u'')$. Then $\tilde{u}'=\iota_{\widetilde{N}}^*(\tilde{u}'')$ and
$\tilde{u}=\tilde{k}^*(\tilde{u}'')$. Since $\tilde{\pi}$ is a
deformation retraction, we have $j\circ\tilde{\pi}\simeq
\iota_{\widetilde{N}}$. Now
\begin{align*}
j_*j_!(a\times b)&=j_*\bigl(\tilde{\pi}_*(\tilde{u}'\cap
\tilde{c}_*(a\times b))\bigr) \\
&={\iota_{\widetilde{N}}}_*
\bigl(\iota_{\widetilde{N}}^*(\tilde{u}'')\cap\tilde{c}_*(a\times b)\bigr) \\
&=\tilde{u}''\cap\tilde{k}_*(a\times b)=\tilde{k}^*(\tilde{u}'')\cap(a\times b) \\
&=\tilde{u}\cap(a\times b).
\end{align*}
This completes the proof. 
\end{proof}

\section{Homotopy Theoretic Description of the Loop Bracket}

Chas and Sullivan constructed the loop bracket on chain level
\cite{CS}. We reformulate their construction in a homotopy theoretic
way, just as Cohen and Jones did in \cite{CJ} for the loop product,
and prove their graded anti-commutation relation.

First we describe the configuration space of two loops on which loop
bracket interaction takes place. We consider the following evaluation
map, where $t\in S^1=\mathbb{R}/\mathbb{Z}$.
\begin{equation}
\begin{gathered}\label{evaluation map for loop bracket}
e: S^1\times LM\times LM \longrightarrow  M\times M, \\
e(t,\gamma,\eta)=
\begin{cases}
\bigl(\gamma(0), \eta(2t)\bigr)& 0\le t\le \tfrac12, \\
\bigl(\gamma(2t-1), \eta(1)\bigr)& \tfrac12\le t\le 1. 
\end{cases}
\end{gathered}
\end{equation}
Using the diagonal map $\phi:M \longrightarrow M\times M$, we consider
a pull-back set $e^{-1}\bigl(\phi(M)\bigr)$ consisting of triples
$(t,\gamma,\eta)$ such that $\gamma(0)=\eta(2t)$ for $0\le
t\le\frac12$, and $\gamma(2t-1)=\eta(1)$ for $\frac12\le t\le 1$.  We
let
\begin{equation}
e^{-1}\bigl(\phi(M)\bigr)_t=
e^{-1}\bigl(\phi(M)\bigr)\cap (\{t\}\times LM\times LM)
\end{equation}
for $0\le t\le 1$. Note that for $t=0,\frac12,1$, the set
$e^{-1}\bigl(\phi(M)\bigr)_t$ describes the same subset
$LM\times_MLM\subset LM\times LM$. Each element of
$e^{-1}\bigl(\phi(M)\bigr)$ is a pair of two loops intersecting at a
point, and are ready to interact to form a single loop. Let
$\iota_1:e^{-1}\bigl(\phi(M)\bigr) \longrightarrow LM$ be an
interaction map given as follows.

\smallskip

\textup{(I)} for $0\le t\le \tfrac12$, $(t,\gamma,\eta)\in e^{-1}\bigl(\phi(M)\bigr)$ if and only if $\gamma(0)=\eta(2t)$ and 
\begin{equation}\label{iota1-1}
\iota_1(t,\gamma,\eta)(s)=
\begin{cases}
\eta(2s)        & 0\le s\le t,\\
\gamma(2s-2t)\quad   & t\le s\le t+\tfrac12,\\
\eta(2s-1)      & t+\tfrac12\le s\le 1.
\end{cases} 
\end{equation}
\textup{(II)} for $\tfrac12\le t\le 1$, $(t,\gamma, \eta)\in e^{-1}\bigl(\phi(M)\bigr)$ if and only if $\gamma(2t-1)=\eta(1)$, and 
\begin{equation}\label{iota1-2}
\iota_1(t,\gamma,\eta)(s)=
\begin{cases}
\gamma(2s)      & 0\le s\le t-\frac12, \\
\eta(2s-2t+1)\quad   & t-\frac12\le s\le t, \\
\gamma(2s-1)    & t\le s\le 1.
\end{cases}
\end{equation}
Thus, for $0\le t\le \frac12$, $\iota_1(t,\gamma,\eta)$ is a loop
starting at the base point of $\eta$, following $\eta$ along its
orientation until it encounters $\gamma$ at $\eta(2t)$, then follow
$\gamma$ from $\gamma(0)$ to $\gamma(1)=\eta(2t)$, then follow along
$\eta$ to $\eta(1)$. Similarly for the case $\frac12\le t \le 1$. This
interaction is exactly the interaction for the loop bracket given in
\cite[\S 4]{CS}. As a function of $t$, during $0\le t\le \frac12$,
$\gamma$ loops move along $\eta$ loops from $\eta(0)$ to $\eta(1)$, then
during $\frac12\le t\le 1$, $\eta$ loops move along $\gamma$ loops
from $\gamma(0)$ to $\gamma(1)$.

Note that although $e^{-1}\bigl(\phi(M)\bigr)_0$ and
$e^{-1}\bigl(\phi(M)\bigr)_{\frac12}$ represent the same set of
configurations of pairs of loops, the effect of $\iota_1$ on these
sets are different. We have $\iota_1(0,\gamma,\eta)=\gamma\cdot\eta$
and $\iota_1(\frac12,\gamma,\eta)=\eta\cdot\gamma$. Note that
$e^{-1}\bigl(\phi(M)\bigr)$ for $t\in[0,\frac12]$ gives a homotopy
between $\gamma\cdot\eta$ and $\eta\cdot\gamma$.

We have the following diagram.
\begin{equation}
\begin{CD}
S^1\times LM\times LM @<{j}<<  e^{-1}\bigl(\phi(M)\bigr) @>{\iota_1}>> LM \\
@V{e}VV   @V{q}VV   @.  \\
M\times M @<{\phi}<< M  @.  
\end{CD}
\end{equation}
where the map $q$ is a restriction of $e$. This diagram defines the
loop bracket.

\begin{defn}
For $a,b\in H_*(LM)$, their loop bracket $\{a,b\}\in
H_{|a|+|b|-d+1}(LM)$ is defined by the following formula
\begin{equation}
\{a,b\}=(-1)^{|a|+d(|a|-d)}
{\iota_1}_*j_!\bigl([S^1]\times a\times b\bigr),
\end{equation}
where $j_!$ is the transfer map associated to the Thom class of $\phi$. 
\end{defn}
The construction of the transfer map $j_!$ is basically the same as
the transfer map $j_!$ appearing in the definition of the loop
product. We go through the construction, and describe aspects
different from the loop product case. Let $\widehat{N}=e^{-1}(N)$ and
let $\hat{c}: S^1\times LM\times LM \longrightarrow
\widehat{N}/\partial\widehat{N}$ be the Thom collapse map. Let
$\hat{u}'\in H^d(\widehat{N},\partial\widehat{N})$ and $\hat{u}\in
H^d(S^1\times LM\times LM)$ be the Thom classes corresponding to Thom
classes $u', u$ of the base manifolds.

We define a lift $\hat{\pi}:\widehat{N} \longrightarrow
e^{-1}\bigl(\phi(M)\bigr)$ of $\pi: N \longrightarrow \phi(M)$
satisfying $q\circ \hat{\pi}=\pi\circ e$ as follows. The construction of the lift $\tilde{\pi}$ can be done abstractly using the homotopy lifting property of the fibration. But here we can be very explicit without difficulty, we give some details.  Let
$(t,\gamma,\eta)\in\widehat{N}$ and let $\pi\circ
e(t,\gamma,\eta)=(z,z)\in\phi(M)$. Let $\ell=(\ell_1,\ell_2):I
\rightarrow M$ be a path in $M\times M$ from $(z,z)$ to
$e(t,\gamma,\eta)\in N$ corresponding to a straight ray from the
origin in a fiber of the vector bundle $\nu$ using $D(\nu)\cong N$. To define
$\hat{\pi}(t,\gamma,\eta)=(t,\gamma_t,\eta_t)\in
e^{-1}\bigl(\phi(M)\bigr)$, we first consider auxiliary loops
$\gamma'_t, \eta'_t$ and modify them.

For $0\le t\le \frac12$, let $\gamma'_t, \eta_t':S^1 \longrightarrow
M$ be loops given as follows
\begin{align*} 
\gamma'_t(s)&=
\begin{cases}
\ell_1(3s)     & 0\le s\le \tfrac13, \\
\gamma(3s-1)   & \tfrac13\le s\le \tfrac23, \\
\ell_1^{-1}(3s-2) & \tfrac23\le s\le 1.
\end{cases}   
\\
\eta'_t(s)&=
\begin{cases}
\eta(3s)   & 0\le s \le \tfrac{2t}3, \\
\ell_2^{-1}(3s-2t)   & \tfrac{2t}3\le s\le \tfrac{2t+1}3, \\
\ell_2(3s-2t-1)    & \tfrac{2t+1}3 \le s\le \tfrac{2t+2}3, \\
\eta(3s-2)     & \tfrac{2t+2}3 \le s\le 1.
\end{cases}
\end{align*}
Here, $\gamma'_t$ is independent of $t$, and $\gamma'_t(0)=z$ for all
$0\le t\le \frac12$. The loop $\eta_t'$ starts at $\eta(0)$ and goes
through $\eta_t'(\frac{2t+1}3)=z$, and comes back to
$\eta'_t(1)=\eta(1)$. To have an element of
$e^{-1}\bigl(\phi(M)\bigr)_t$, we need to rotate $\eta'_t$. We let
$\gamma_t(s)=\gamma'_t(s)$ and
$\eta_t(s)=\eta'_t(s-\frac{4t-1}3)$. Then
$\eta_t(2t)=z=\gamma_t(0)$. Thus, $(t,\gamma_t,\eta_t)\in
e^{-1}\bigl(\phi(M)\bigr)_t$, and we define
$\hat{\pi}(t,\gamma,\eta)=(t,\gamma_t,\eta_t)$ for $0\le t\le
\frac12$.

For $\frac12\le t\le 1$, define $\gamma'_t, \eta'_t : S^1
\longrightarrow M$ as follows.
\begin{align*}
\gamma'_t(s)&=
\begin{cases}
\gamma(3s)     & 0\le s\le \tfrac{2t-1}3, \\
\ell_1^{-1}(3s-2t+1)    & \tfrac{2t-1}3\le s\le \tfrac{2t}3, \\
\ell_1(3s-2t)   & \tfrac{2t}3 \le s\le \tfrac{2t+1}3, \\
\gamma(3s-2)   & \tfrac{2t+1}3 \le s\le 1.
\end{cases}
\\
\eta'_t(s)&=
\begin{cases}
\ell_2(3s)   & 0\le s\le \tfrac13, \\
\eta(3s-1)  & \frac13\le s\le \tfrac23, \\
\ell_2^{-1}(3s-2)  & \tfrac23 \le s \le 1.
\end{cases}
\end{align*} 
The loop $\gamma'_t$ starts at $\gamma(0)$, goes through
$\gamma'_t(\frac{2t}3)=z$, and comes back to $\gamma(1)$. Thus, we let
$\gamma_t(s)=\gamma'_t(s-\frac{4t-3}3)$ by rotating $\gamma'_t$, and
we let $\eta_t=\eta'_t$ for $\frac12\le t\le 1$. We then have
$\gamma_t(2t-1)=z=\eta_t(0)$, and $(t,\gamma_t,\eta_t)\in
e^{-1}\bigl(\phi(M)\bigr)_t$. Hence we define
$\hat{\pi}(t,\gamma,\eta)=(t,\gamma_t,\eta_t)$ for $\frac12\le t\le
1$. Since at $t=0,\frac12, 1$, we have
$(\gamma_t,\eta_t)=(\ell_1\gamma\ell_1^{-1}, \ell_2\eta\ell_2^{-1})$,
these two families paste together to define a map $\hat{\pi}:
\widehat{N} \longrightarrow e^{-1}\bigl(\phi(M)\bigr)$. By considering
partial path $\ell_{[t,1]}$, we see that $\hat{\pi}$ is a deformation
retraction.

We define the transfer map $j_!$ as the composition of the following maps:
\begin{equation}
j_!: H_*(S^1\times LM\times LM) \xrightarrow{\hat{c}_*}
H_*(\widehat{N},\partial\widehat{N}) \xrightarrow{\hat{u}'\cap(\ )}
H_{*-d}(\widehat{N}) \xrightarrow[\cong]{\hat{\pi}_*}
H_{*-d}\bigl(e^{-1}(\phi(M))\bigr).
\end{equation}
Geometrically, $j_!$ corresponds to taking the transversal
intersection of a cycle representing $[S^1]\times a\times b\in
H_*(S^1\times LM\times LM)$ with $e^{-1}\bigl(\phi(M)\bigr)$, and
consider its homology class with appropriate orientation in
$H_{*}\bigl(e^{-1}(\phi(M))\bigr)$.

Next, we give a homotopy theoretic proof of the graded
anti-commutation relation for the loop bracket.

\begin{prop}
For $a,b\in H_*(LM)$, the anti-commutation relation for the loop
bracket is given in the following form.
\begin{equation}
\{a,b\}=-(-1)^{(|a|-d+1)(|b|-d+1)}\{b,a\}.
\end{equation}
\end{prop} 
\begin{proof} We have the following commutative diagram whose 
commutativity can be directly checked from definition.
\begin{equation*}
\begin{CD}
M\times M @<{e}<< S^1\times LM\times LM @<{j}<< e^{-1}\bigl(\phi(M)\bigr)  
@>{\iota_1}>>   LM   \\
@V{T}VV   @V{R_{\frac12}\times T}VV    @V{R_{\frac12}\times T}VV @|  \\
M\times M @<{e}<<   S^1\times LM\times LM @<{j}<< e^{-1}\bigl(\phi(M)\bigr)
@>{\iota_1}>>  LM 
\end{CD}
\end{equation*}
Here $R_{\frac12}$ is a rotation of loops by $\frac12$, and we let
$\widehat{T}=R_{\frac12}\times T$. In the associated homology square
with transfers $j_!$, the middle square commutes up to a sign. To determine
the sign, we compare $j_*j_!\widehat{T}_*$ and
$j_*\widehat{T}_*j_!$. On the one hand, using Proposition
\ref{jj-formula} we have $j_*j_!\widehat{T}_*([S^1]\times a\times
b)=\hat{u}\cap \widehat{T}_*([S^1]\times a\times b)$. On the other
hand, again using Proposition \ref{jj-formula},
\begin{align*}
j_*\widehat{T}_*j_!([S^1]\times a\times b)
&= \widehat{T}_*j_*j_!([S^1]\times a\times b)\\
&=\widehat{T}_*\bigl(\hat{u}\cap ([S^1]\times a\times b)\bigr)\\
&=\widehat{T}^*(\hat{u})\cap\widehat{T}_*([S^1]\times a\times b).
\end{align*}
Since the left square of the diagram commutes, and the Thom class $u$
satisfies $T^*(u)=(-1)^du$, we have
$\widehat{T}^*(\hat{u})=(-1)^d\hat{u}$. Hence
$j_*\widehat{T}_*j_!=(-1)^dj_*j_!\widehat{T}_*$, consequently,
$\widehat{T}_*j_!=(-1)^dj_!\widehat{T}_*$. Now,
\begin{align*}
\{a,b\}&= (-1)^{|a|+d(|a|-d)}{\iota_1}_*j_!([S^1]\times a\times b)  \\
&=(-1)^{|a|+d(|a|-d)+d}{\iota_1}_*j_!\widehat{T}_*([S^1]\times a\times b)\\
&=(-1)^{|a||b|+|a|+d|a|}{\iota_1}_*j_!([S^1]\times b\times a)  \\
&=-(-1)^{(|a|-d+1)(|b|-d+1)}\{b,a\}.
\end{align*}
This completes the proof of the anti-commutativity of the loop bracket. 
\end{proof}

\section{BV-operator and Derivation}

We examine interaction diagrams corresponding to operations which
assign $a\cdot\Delta b$ and $(\Delta a)\cdot b$ to $a,b\in
H_*(LM)$. The relevant diagrams are
\begin{gather}
S^1\times LM\times LM \xrightarrow{T\times 1} LM\times S^1\times LM
\xrightarrow{1\times\Delta} LM\times LM \xleftarrow{j} LM\times_M LM
\xrightarrow{\iota} LM, \label{(i)} \\ 
S^1\times LM\times LM
\xrightarrow{\Delta\times 1} LM\times LM \xleftarrow{j} LM\times_M LM
\xrightarrow{\iota} LM. \label{(ii)}
\end{gather}
For $a,b\in H_*(LM)$, these diagrams give 
\begin{align}
\iota_*j_!(1\times\Delta)_*(T\times 1)_*([S^1]\times a\times b)&=
(-1)^{|a|+d(|a|-d)}a\cdot\Delta b, \label{a-Delta b} \\
\iota_*j_!(\Delta\times 1)_*([S^1]\times a\times b)&=
(-1)^{d|a|}(\Delta a)\cdot b.
\end{align}
The diagram \eqref{(i)} fits into the following commutative diagram:
\begin{equation*}
\begin{CD} 
S^1\times LM\times LM @<{j'}<<   e_1^{-1}\bigl(\phi(M)\bigr)   
@>{\iota'}>> LM \\
@V{T\times 1}VV   @V{T\times 1}VV  @|  \\
LM\times S^1\times LM @<{j'''}<<  LM\times_M(S^1\times LM)   @>>> LM  \\
@V{1\times \Delta}VV  @VVV  @|  \\
LM\times LM  @<{j}<<  LM\times_M LM @>{\iota}>>   LM  \\
@V{p\times p}VV  @V{p}VV   @. \\
M\times M @<{\phi}<< M @.  
\end{CD}
\end{equation*} 
where $e_1=(p\times p)(1\times \Delta)(T\times 1)$ and is given by
$e_1(t,\gamma,\eta)=(\gamma(0), \eta(t))$, and 
\begin{equation}
\iota'(t,\gamma,\eta)(s)=
\begin{cases} 
\gamma(2s)  & 0\le s \le \tfrac12, \\
\eta(2s-1+t)   & \tfrac12 \le s\le 1.
\end{cases} 
\end{equation}
The loop $\iota'(t,\gamma,\eta)$ starts at $\gamma(0)$ and follows the
orientation of $\gamma$ all the way to $\gamma(1)=\eta(t)$, then follows 
the entire $\eta$ from $\eta(t)$ to $\eta(t+1)$. 

For maps $j, j', j'''$, by using the pull-backs of the same Thom class
$u\in H^d(M\times M)$, the resulting transfer maps $j_!, j'_!, j_!'''$
are all compatible and the induced homology diagram with these
transfer maps commutes. In particular, \eqref{a-Delta b} gives 
\begin{equation}\label{derivation 1}
\iota'_*j_!'([S^1]\times a\times b)
=(-1)^{|a|+d(|a|-d)} a\cdot\Delta b.
\end{equation} 

The diagram \eqref{(ii)} fits into the following commutative diagram.
\begin{equation}
\begin{CD}
S^1\times LM\times LM @<{j''}<<  e_2^{-1}\bigl(\phi(M)\bigr)  @>{\iota''}>>  LM \\
@V{\Delta\times 1}VV   @VVV    @|   \\
LM\times LM  @<{j}<< LM\times_M LM  @>{\iota}>>  LM \\
@V{p\times p}VV  @VVV   @.   \\
M\times M  @<{\phi}<<  M @.  
\end{CD}
\end{equation}
where $e_2=(p\times p)(\Delta\times 1)$ is given by 
$e_2(t,\gamma, \eta)=\bigl(\gamma(t),\eta(0)\bigr)$, and 
$e_2^{-1}\bigl(\phi(M)\bigr)$ consists of $(t,\gamma,\eta)$ such that $\gamma(t)=\eta(0)$. Then the map $\iota''$ is given by 
\begin{equation}
\iota''(t,\gamma,\eta)(s)=
\begin{cases} 
\gamma(2s+t)  &  0\le s\le \tfrac12, \\
\eta(2s-1)   & \tfrac12\le s\le 1.
\end{cases}
\end{equation}
The loop $\iota''(t,\gamma,\eta)$ starts at $\gamma(t)$ and follows the orientation of $\gamma$ to $\gamma(t+1)=\eta(0)$, and then moves along $\eta$ from $\eta(0)$ to $\eta(1)=\gamma(t)$. 

Transfer maps $j_!,j_!''$ can be constructed using pull-backs of the same Thom class $u\in H^d(M\times M)$. Then the induced homology diagram with transfers commutes, and we have 
\begin{equation}\label{derivation 2}
\iota''_*j''_!\bigl([S^1]\times a\times b\bigr)
=(-1)^{d|a|}(\Delta a)\cdot b.
\end{equation} 
To construct the loop bracket, we used the evaluation map $e: S^1\times LM\times LM \longrightarrow M\times M$ given in \eqref{evaluation map for loop bracket}. Now we note that the evaluation maps $e_1$ and $e_2$ are precisely the first half and the second half of $e$. Namely, 
\begin{equation}
e(t,\gamma,\eta)=
\begin{cases}
e_1(2t,\gamma,\eta) &  0\le t\le \tfrac12,\\
e_2(2t-1,\gamma,\eta)    & \tfrac12 \le t\le 1.
\end{cases}
\end{equation}
Thus, we combine $\iota'$ and $\iota''$ to define $\iota_2$ by 
\begin{equation}
\begin{gathered}
\iota_2: e^{-1}\bigl(\phi(M)\bigr) \longrightarrow LM  \\
\iota_2(t,\gamma,\eta)=
\begin{cases}
\iota'(2t,\gamma,\eta)  & 0\le t\le \tfrac12, \\
\iota''(2t-1,\gamma,\eta)  &  \tfrac12 \le t\le 1.
\end{cases}
\end{gathered}
\end{equation}
Thus the following diagram combines \eqref{(i)} and \eqref{(ii)}. 
\begin{equation}
S^1\times LM\times LM \xleftarrow{j}  e^{-1}\bigl(\phi(M)\bigr) 
\xrightarrow{\iota_2} LM .
\end{equation}
This diagram gives what we expect. 
\begin{prop} 
For $a,b\in H_*(LM)$, we have 
\begin{equation}
{\iota_2}_*j_!([S^1]\times a\times b)=(-1)^{|a|+d(|a|-d)}a\cdot\Delta b +
(-1)^{d|a|}\Delta a\cdot b.
\end{equation}
\end{prop} 
\begin{proof} We introduce some notations. Let $I_1=[0,\frac12]$, $I_2=[\frac12,1]$, $S_1^1=I_1\partial I_1$, $S^1_2=I_2\partial I_2$, and $q: S^1 \rightarrow S^1_1\vee S^1_2$ be an identification map.

Since $e(0,\gamma,\eta)=e(\frac12,\gamma,\eta)=(\gamma(0),\eta(0))\in M\times M$,  
The map $e$ factors through $(S^1_1\vee S^1_2)\times LM\times LM$. We consider the following diagram. 
\begin{equation*}
\begin{CD}
S^1\times LM\times LM  @<{j'}<<  e_1^{-1}\bigl(\phi(M)\bigr)
@>{\iota'}>>   LM  \\
@V{r_1\times 1\times 1}VV   @V{r_1'}VV  @|  \\  
(S^1_1\vee S^1_2)\times LM\times LM   @<{\hat{j}}<< e^{-1}\bigl(\phi(M)\bigr)
@>{\iota_2}>>  LM  \\
@A{r_2\times 1\times 1}AA  @A{r_2'}AA   @|  \\
S^1\times LM\times LM  @<{j''}<<  e_2^{-1}\bigl(\phi(M)\bigr)
@>{\iota''}>>   LM
\end{CD}
\end{equation*}
where $r_i: S^1 \rightarrow S^1_i$ for $i=1,2$ are given by $r_1(t)=\frac{r}2$, $r_2(t)=\frac{t+1}2$, and $r_i'$ for $i=1,2$ are restrictions of $r_i\times 1\times 1$. Since $\hat{j}=(q\times 1\times 1)j$, using pull-backs of the same Thom class $u$ from $M\times M$, we have $j_!=\hat{j}_!(q\times 1\times 1)_*$. Similarly, the homology diagram with transfers $j_!, \hat{j}_!, j''_!$ induced from the above diagram commutes. Hence 
\begin{align*}
{\iota_2}_*j_!([S^1]\times a\times b)&=
{\iota_2}_*\hat{j}_!([S^1_1]\times a\times b+[S^1_2]\times a\times b) \\
&=\iota'_*j'_!([S^1]\times a\times b) 
+\iota''_*j''_!([S^1]\times a\times b) \\
&=(-1)^{|a|+d(|a|-d)}a\cdot\Delta b
+(-1)^{d|a|}\Delta a\cdot b,
\end{align*}
using \eqref{derivation 1} and \eqref{derivation 2}.
This completes the proof. 
\end{proof} 

For convenience, we write out the map $\iota_2$ explicitly. 

(I) For $0\le t\le \tfrac12$, $(t,\gamma,\eta)\in e^{-1}\bigl(\phi(M)\bigr)$ if and only if $\gamma(0)=\eta(2t)$ and 
\begin{equation}\label{iota2-1}
\iota_2(t,\gamma,\eta)(s)=
\begin{cases}
\gamma(2s)    &   0\le s\le \tfrac12, \\
\eta(2s-1+2t)   & \tfrac12 \le s\le 1. 
\end{cases}
\end{equation}

(II) For $\frac12\le t\le 1$, $(t,\gamma,\eta)\in e^{-1}\bigl(\phi(M)\bigr)$
if and only if $\gamma(2t-1)=\eta(1)$ and 
\begin{equation}\label{iota2-2}
\iota_2(t,\gamma,\eta)(s)=
\begin{cases}
\gamma(2s+2t-1)   &  0\le s\le \tfrac12, \\
\eta(2s-1)   &   \tfrac12 \le s \le 1.
\end{cases}
\end{equation}

\section{A Proof of the BV Identity}

We combine the descriptions of $\{a, b\}$ and $\Delta a\cdot b+(-1)^{|a|-d}a\cdot\Delta b$ in previous sections to prove the BV identity  
\begin{equation}
\Delta(a\cdot b)=\Delta a\cdot b + (-1)^{|a|-d}a\cdot \Delta b 
-(-1)^{|a|-d}\{a.b\}. 
\end{equation} 
The minus sign in front of the loop bracket is due to our choice of $S^1$ action $\Delta: S^1\times LM \longrightarrow LM$ given by $\Delta(t,\gamma)=\gamma_t$, where $\gamma_t(s)=\gamma(s+t)$. If we use the opposite action $\Delta'$ given by $\Delta'(t,\gamma)=\gamma_{(-t)}$, then with respect to the action, we get the plus sign in front of the loop bracket in the above BV identity. 

Results in previous sections can be summarized by the following diagram and identities for $a,b\in H_*(LM)$: 
\begin{gather}
S^1\times LM\times LM \xleftarrow{j} e^{-1}\bigl(\phi(M)\bigr) 
\xrightarrow{\iota_1,\iota_2}   LM    \\
{\iota_1}_*j_!([S^1]\times a\times b)=(-1)^{|a|+d(|a|-d)}\{a,b\}, \label{loop bracket}\\
{\iota_2}_*j_!([S^1]\times a\times b)=(-1)^{|a|+d(|a|-d)}a\cdot\Delta b
+(-1)^{d|a|}\Delta a\cdot b.\label{bv derivation}
\end{gather}
Note that the above two interactions are defined on the same configuration set $e^{-1}\bigl(\phi(M)\bigr)$, and the only difference between the loop bracket and the BV derivation lies in the difference of $\iota_1$ and $\iota_2$, which turns out to be a simple reparametrization of loops. To describe this, let $\pi:e^{-1}\bigl(\phi(M)\bigr) \xrightarrow{j} S^1\times LM\times LM  \xrightarrow{\pi_1} S^1$ be the projection map onto the $S^1$ factor. 
\begin{lem}\label{rotation}
Let $\widehat{\Delta} : S^1\times LM \longrightarrow S^1\times LM $ be given by $\widehat{\Delta}(t,\gamma)=(t,\gamma_t)$. Then the following diagram commutes. 
\begin{equation}\label{iota1-iota2}
\begin{CD}
e^{-1}\bigl(\phi(M)\bigr) @>{(\pi,\iota_1)}>> S^1\times LM \\
@|  @V{\widehat{\Delta}}VV  \\
e^{-1}\bigl(\phi(M)\bigr) @>{(\pi,\iota_2)}>> S^1\times LM 
\end{CD}
\end{equation}
In other words, $\iota_1(t,\gamma,\eta)_t=\iota_2(t,\gamma,\eta)$, for $(t,\gamma,\eta)\in e^{-1}\bigl(\phi(M)\bigr)$. 
\end{lem} 
\begin{proof}  This is straightforward checking using \eqref{iota1-1}, \eqref{iota1-2}, \eqref{iota2-1}, and \eqref{iota2-2}. 

When $0\le t\le\frac12$, we have 
\begin{equation*}
\iota_1(t,\gamma,\eta)_t(s)=\iota_1(t,\gamma,\eta)(s+t)=
\left\{
\begin{aligned}
&\gamma(2s)  & &0\le s\le \tfrac12 \\
&\eta(2s+2t-1)  & &\tfrac12\le s\le 1
\end{aligned}
\right\}
=\iota_2(t,\gamma,\eta)(s).
\end{equation*}

When $\frac12\le t\le 1$, we have 
\begin{equation*}
\iota_1(t,\gamma,\eta)_t(s)=\iota_1(t,\gamma,\eta)(s+t)=
\left\{
\begin{aligned}
&\gamma(2s+2t-1) & & 0\le s\le \tfrac12 \\
&\eta(2s-1)  & &\tfrac12\le s\le 1
\end{aligned}
\right\}
=\iota_2(t,\gamma,\eta)(s).
\end{equation*}
This completes the proof.
\end{proof} 

To study the homological behavior of the above diagram, we need to know the Thom class for the embedding $j$.
\begin{prop} 
The Thom class $\hat{u}=e^*(u)$ of the embedding $j$ in the following diagram 
\begin{equation}
\begin{CD}
S^1\times LM\times LM  @<{j}<<  e^{-1}\bigl(\phi(M)\bigr)  \\
@V{e}VV   @VVV  \\
M\times M @<{\phi}<<  M 
\end{CD}
\end{equation}
is given by 
\begin{equation}\label{eu}
\hat{u}=e^*(u)=\Delta^*(\tilde{u})=1\times \tilde{u} + \{S^1\}\times \Delta\tilde{u}, 
\end{equation}
where $\Delta: S^1\times L(M\times M) \longrightarrow L(M\times M)$ is the $S^1$ action on the free loop space $L(M\times M)$, and $\tilde{u}=(p\times p)^*(u)$ is the pull-back of the Thom class $u$ of the diagonal map $\phi$ to $LM\times LM$.  
\end{prop} 
\begin{proof} We consider two loops $\delta_i: S^1 \rightarrow S^1\times S^1$ for $i=1,2$ given by 
\begin{equation*}
\delta_1(t)=
\begin{cases}
(0,2t)  & 0\le t\le \tfrac12 \\
(2t-1,1)  & \tfrac12 \le t\le 1
\end{cases},
\quad \delta_2(t)=(t,t), 0\le t\le 1.
\end{equation*}
The loop $\delta_1$ goes around he first circle, then around 
the second circle, and the loop $\delta_2$ is the diagonal loop. 
Obviously, these two loops are homotopic to each other. 
We consider the following composition map for $i=1,2$:   
\begin{multline} 
S^1\times LM\times LM \xrightarrow{\delta_i\times 1\times 1}  S^1\times S^1\times LM\times LM \\
\xrightarrow{1\times T\times 1} S^1\times LM\times S^1\times LM 
\xrightarrow{\Delta\times \Delta}  LM\times LM \xrightarrow{p\times p}  M\times M.
\end{multline}
Then for $i=1$, the above map is exactly the evaluation map $e$. For $i=2$, the above composition is the same as $S^1\times L(M\times M) \xrightarrow{\Delta} L(M\times M) \xrightarrow{p} M\times M$. Since $\delta_1$ and $\delta_2$ are homotopic to each other, we have $e^*(u)=\Delta^*p^*(u)=\Delta^*(\tilde{u})$. This completes the proof. 
\end{proof} 
\begin{prop}\label{Pi Iota1}
In the diagram 
\begin{equation}
S^1\times LM\times LM  \xleftarrow{j} e^{-1}\bigl(\phi(M)\bigr) 
\xrightarrow{(\pi,\iota_1)}  S^1\times LM,
\end{equation}
for $a,b\in H_*(LM)$ homological behavior is given by 
\begin{equation}\label{pi-iota1}
(\pi,\iota_1)_*j_!([S^1]\times a\times b)=[S^1]\times (-1)^{d|a|}a\cdot b 
+[0]\times (-1)^{|a|+d(|a|-d)}\{a,b\},
\end{equation}
where $0\in S^1$ is the base point of $S^1$. 
\end{prop} 
\begin{proof}  Let $(\pi,\iota_1)_*j_!([S^1]\times a\times b)=[S^1]\times x + 
[0]\times y$ for some $x,y\in H_*(LM)$. If $\pi_2:S^1\times LM \rightarrow LM$ is the projection onto the second factor, then we have 
\begin{equation*} 
y={\pi_2}_*(\pi,\iota_1)_*j_!([S^1]\times a\times b)={\iota_1}_*j_!([S^1]\times a\times b)=(-1)^{|a|+d(|a|-d)}\{a,b\}.
\end{equation*} 
To identify $x$, let $h:\{0\} \rightarrow S^1$ be the inclusion map, and consider the following homology diagram.
\begin{equation*}
\begin{CD}
H_*(S^1\times LM^2) @>{j_!}>> H_{*-d}\bigl(e^{-1}\bigl(\phi(M)\bigr)\bigr)
@>{(\pi,\iota_1)_*}>>  H_{*-d}(S^1\times LM)  \\
@V{(h\times 1\times 1)_!}VV   @V{(h\times 1\times 1)_!}VV  @V{(h\times 1)_!}VV \\
H_{*-1}(\{0\}\times LM^2)  @>{(1\times j)_!}>>  
H_{*-d-1}(\{0\}\times LM\underset{M}{\times} LM)  @>{(1\times\iota)_*}>> 
H_{*-d-1}(\{0\}\times LM) 
\end{CD}
\end{equation*} 
We show that the left homology square commutes up to $(-1)^d$. Since the diagram commutes on space level, the homology diagram with transfers commutes up to a sign. We determine this sign. Since $e\circ(h\times 1\times 1)=p\times p: LM\times LM \rightarrow M\times M$, we have $(h\times 1\times 1)^*e^*(u)=(p\times p)^*(u)=\tilde{u}$. Also, since the Thom class of the inclusion $h:\{0\} \rightarrow
S^1$ is given by $\{S^1\}$, we have 
\begin{equation*}
(h\times 1\times 1)_*(h\times 1\times 1)_!([S^1]\times a\times b)
=(\{S^1\}\times 1\times 1)\cap([S^1]\times a\times b)
\end{equation*} 
on $H_*(S^1\times LM\times LM)$.  Now 
\begin{align*}
(h\times 1\times 1)_*&(1\times j)_*(1\times j)_!(h\times 1\times 1)_!
([S^1]\times a\times b)  \\
&=(h\times 1\times 1)_*\bigl(\tilde{u}\cap(h\times 1\times 1)_!
([S^1]\times a\times b)\bigr)  \\
&=e^*(u)\cap\bigl((\{S^1\}\times1\times 1)\cap([S^1]\times a\times b)\bigr) \\
&=e^*(u)\cap([0]\times a\times b)=[0]\times \bigl(\tilde{u}\cap(a\times b)\bigr), 
\end{align*}
where in the last identity, we used \eqref{eu}. On the other hand, 
\begin{align*}
(h\times 1\times 1)_*&(1\times j)_*(h\times 1\times1)_!j_!([S^1]\times a\times b)\\
&=j_*(h\times 1\times 1)_*(h\times 1\times 1)_!j_!([S^1]\times a\times b)  \\
&=j_*\bigl(j^*(\{S^1\}\times 1\times 1)\cap j_!([S^1]\times a\times b)\bigr)  \\
&=(\{S^1\}\times 1\times 1)\cap j_*j_!([S^1]\times a\times b)  \\
&=(\{S^1\}\times 1\times 1)\cap \bigl(e^*(u)\cap ([S^1]\times a\times b)\bigr)
\end{align*}
Since $e^*(u)=1\times\tilde{u} + \{S^1\}\times\Delta\tilde{u}$ by \eqref{eu}, we have 
$(\{S^1\}\times 1\times 1)\cup e^*(u)=\{S^1\}\times\tilde{u}$. Hence the last formula above is equal to $(-1)^d[0]\times \bigl(\tilde{u}\cap(a\times b)\bigr)$. Thus comparing the above two computations, we have $(h\times 1\times 1)_!j_!=(-1)^d(1\times j)_!(h\times 1\times 1)_!$, and the left square of the homology diagram commutes up to $(-1)^d$. Thus, the homology diagram implies 
\begin{align*}
[0]\times x&=(h\times 1)_!(\pi,\iota_1)_*j_!([S^1]\times a\times b) \\
&=(-1)^d[0]\times \iota_*j_!(a\times b)=[0]\times (-1)^{d|a|}a\cdot b.
\end{align*}
Hence $x=(-1)^{d|a|} a\cdot b$. This completes the proof. 
\end{proof} 

\begin{thm}\label{proof of BV}  Let $\Delta: S^1\times LM \rightarrow LM$ be the $S^1$ action map given by $\Delta(t,\gamma)=\gamma_t$, where $\gamma_t(s)=\gamma(s+t)$ for $s,t\in S^1=\mathbb{R}/\mathbb{Z}$. Then for $a,b\in H_*(LM)$, the BV identity holds. 
\begin{equation}
\Delta(a\cdot b)=(\Delta a)\cdot b + (-1)^{|a|-d}a\cdot\Delta b -(-1)^{|a|-d}\{a,b\}. 
\end{equation}
\end{thm}
\begin{proof} By \eqref{iota1-iota2}, we have $\iota_2=\Delta\circ(\pi,\iota_1)$. Hence identities \eqref{bv derivation} and \eqref{pi-iota1} imply 
\begin{align*}
(-1)^{|a|+d(|a|-d)}a\cdot &\Delta b + (-1)^{d|a|}(\Delta a)\cdot b
={\iota_2}_*j_!([S^1]\times a\times b) \\
&={\Delta}_*(\pi,\iota_1)_*j_!([S^1]\times a\times b) \\
&={\Delta}_*\bigl([S^1]\times (-1)^{d|a|}a\cdot b + 
[0]\times (-1)^{|a|+d(|a|-d)}\{a,b\}\bigr) \\
&=(-1)^{d|a|}\Delta(a\cdot b) + (-1)^{|a|+d(|a|-d)}\{a,b\}.
\end{align*}
Hence canceling some signs, we get 
\begin{equation*}
\Delta(a\cdot b) + (-1)^{|a|-d}\{a,b\} 
=(\Delta a)\cdot b + (-1)^{|a|-d}a\cdot \Delta b.
\end{equation*}
This completes the proof of BV identity. 
\end{proof}

\end{document}